\documentclass[Svgc]{Svgc}
\usepackage{graphics}
\usepackage{amsfonts}

\begin{document}
\title{Isotemporal classes of diasters, beachballs, and daisies}
\author{Benjamin de Bivort }                     
\titlerunning{Isotemporal diasters, beachballs, and daisies}%
\authorrunning{Benjamin de Bivort}%
\institute{Harvard University, Department of Organismic and Evolutionary Biology 
\\
Northwest Building Room 235.30
\\
52 Oxford St
\\
Cambridge, MA 02138, USA
\\
debivort@oeb.harvard.edu}
\maketitle

\begin{abstract} If the vertices composing a network interact at distinct time points, the temporal ordering of these
interactions and the network's graph structure are sufficient to convey the routes by which information can flow in the
network. Two networks with real-valued edge labels are \emph{temporally isomorphic} if there exists  a graph isomorphism
$\phi:N\to M$ that preserves \emph{temporal paths} --- paths in which sequential edge labels are strictly increasing. An
equivalence class of temporally ismorphic networks is known as an \emph{isotemporal class}. Methods to determine the
number of isotemporal classes of a particular graph structure ($\mathcal{N}(G)$) are non-obvious, and refractory to
traditional techniques such as P\'{o}lya enumeration, P\'{o}lya (1937). Here, I present a simple formula for the number
of isotemporal classes of \emph{diasters}, graphs composed of a vertex of degree $a+1$ connected to a vertex of degree
$b+1$, with all other vertices of degree 1 (denoted $D(a,b)$). In particular, $\mathcal{N}(D(a,b)))=ab+a+b+1$ if $a\neq
b$, and $\mathcal{N}(D(a,a)))=\frac{1}{2} (a^2+3a+2)$ otherwise. This formula is then extended to five additional types
of pseudograph by application of a theorem that states $\mathcal{N}(G)$ is preserved between two graph types if edge
adjacencies and automorphisms are preserved, and provided that any two networks are members of the same isotemporal
class if and only if they are isomorphic by transpositions of sequential edge labels on non-adjacent edges.
\end{abstract}

\section{Introduction}

Metaphorically interpreting numerical labels on graph edges as the time of an interaction between the vertices is a
framework that can be used to model a number of real systems. For example, the history of European 
warfare could be recorded on a temporal network in which vertices correspond to nations, edges conflicts, and 
edge labels the year of the combat. Bearman \emph{et al.} (2004) analyzed the romantic interactions of 800 American
adolescents, creating a sexual exploit temporal network with 800 vertices and 477 edges, each labeled with the date of
the relationship. Despite the clearly utility of using edge labels as times, analytic treatment of this topic is
restricted to a small number of results such as Kempe \emph{et al.} (2002). 

The framework of a temporal network leads to a plethora of intriguing theoretical questions. For example, if a 
signal or object (e.g. a virus within students) can propagate through the network by crossing an edge only at its 
labeled time, what kind of temporal connectivity is possible in networks of a particular structure, or random 
networks? Moreover, what variety of edge interaction orderings have the same temporal connectivity, and how 
many different types of ``fundamental temporal structures'' (isotemporal classes) exist (where two differently 
labeled networks are considered to have the same structure if they have isomorphic temporal connectivity)? 

In a parallel study (de Bivort (2007)), closed formulas for the number of isotemporal classes of $n$-gons 
($n$-cycles) are determined. These formulas depend heavily on the number and types of symmetries for each $n$ and
contain several terms superficially similar to the closed formula for the number of rotationally-distinct 2-color
necklaces, Fine (1958). In this case, the type of graph under consideration is naturally characterizable by a single
parameter, $n$ (the length of the cycle). In this paper, I report closed formulas for the number of 
isotemporal classes of a set of graphs best characterized by two parameters, and then show that these results can 
be extended to at least five other two-parameter types of pseudographs, using mappings that preserve the number 
isotemporal classes, but not necessarily the vertex structure of the initial graph.

\section{Definitions}

\begin{definition} --- a \textbf{temporal network} $N=\{V,E,T,\tau\}$ is a collection of vertices $V$, 
edges $E = \{\{v_i,v_j\} | v_i, v_j \in V \}$, temporal labels $T \subset \mathbb{R}$ with $|T|\leq|E|$, 
and a temporal mapping $\tau : E \to T$.
\end{definition}

For the rest of this work, we will make the assumption that $\tau$ is a bijection, i.e. $|T|=|E|$. Networks 
of this sort can be considered ``serial temporal networks'' because every interaction happens at a distinct time, and
therefore the edge events happen strictly in series, as opposed to allowing $|T|<|E|$ which generates a ``parallel
temporal network'' in which multiple events happen simultaneously. For convenience we will also assume for convenience that
the temporal values composing $T$ are $\{1,2, \ldots, t\}$ where $t$ is the total number of edges.

The details of a network's temporal connectivity can be considered by enumerating the various temporal paths through the
network:
\begin{definition}
\label{TemporalPathDefinition} --- an ordered series of edges $\langle \{v_1,v_2\},\{v_2,v_3\}, \ldots,$
$\{v_{k},v_{k+1}\}\rangle$ where $\tau(\{v_a,v_{a+1}\}) < \tau(\{v_b,v_{b+1}\})$ for all $a<b$, is a \textbf{temporal
path} of length $k$.
\end{definition}

The set of a temporal network's temporal paths does not uniquely determine its labeling $\tau$. Indeed it is the
exceptional case when there is a single labeling with a particular set of temporal paths. Figure 1 shows two networks
whose temporal paths are preserved by switching temporal labels $1$ and $2$. Each pair of networks in the example are
non superimposable when edge labels are considered, and therefore have different temporal mappings. However, because all
temporal paths are preserved, they are temporally isomorphic.

\begin{definition} --- two temporal networks $N=\{V,E,T,\tau\}$ and $M=\{V',E',$ $T,\tau'\}$ are \textbf{temporally
isomorphic}, denoted $N\cong_\mathrm{T}M$), if there exists a mapping $\phi:V \to V'$ such that 
\begin{enumerate}
\item if $\{v_1,v_2\} \in E$ then $\{\phi(v_1),\phi(v_2)\} \in E'$ and 
\item if $P=\langle \{v_1,v_2\},$ $\{v_2,v_3\}, \ldots,$ $\{v_{k},v_{k+1}\}\rangle$ is a temporal path in $N$ then
$\langle \{\phi(v_1),\phi(v_2)\},$ $\{\phi(v_2),\phi(v_3)\}, \ldots,$ $\{\phi(v_{k}),\phi(v_{k+1})\}\rangle$ is a
temporal path in $M$.
\end{enumerate}
\end{definition}

If $M$ and $N$ are at least graphically isomorphic, we will denote this as $N\cong_\mathrm{G}M$. If $\phi:M\to N$ is a
graphical isomorphism that preserves temporal labels --- for all edges, $\tau(\{v_i,v_j\})=\tau'(\{ \phi(v_i), \phi(v_j)
\})$ --- then $M$ and $N$ are \textbf{label} isomorphic ($N\cong_\mathrm{L}M$).

\begin{definition} --- The set of all $N$ such that $N\cong_\mathrm{T}M$ for a particular $M$ is an
\textbf{isotemporal class}.
\end{definition}

For $N=\{V,E,T,\tau\}$, we will let $\mathcal{N}(N)$ indicate the total number of isotemporal classes of networks
graphically isomorphic to $N$ for all $\tau$. 

\section{Isotemporal classes of diasters}

\begin{definition} --- A \textbf{diaster} is any graph that contains two connected vertices (the
\textbf{central edge}) and any number of vertices of degree one attached to those two vertices (
\textbf{peripheral edges}).
\end{definition}

The class of diasters is best described by two parameters --- the number of degree-1 vertices attached to one of the
internal vertices ($a$), and the number attached to the other ($b$). This is denoted $D(a,b)=D(b,a)$. All diasters
smaller than $D(3,4)$ are shown in Figure 2. Diasters are trees with at most two internal branching points and
an arbitrary number of leaves. 

The number of isotemporal classes of an arbitrary diaster ($\mathcal{N}(D(a,b))$) is highly restricted by its structure.
The longest temporal path in any diaster will contain at most three edges, since the requirement that subsequent edges
of a temporal path have temporal values greater than the current edge excludes the possibility of doubling back across
an edge that has already been crossed. Moreover, because the central edge is necessarily the second edge in any path of
length greater than one, to enumerate all the isotemporal classes of a diaster it is sufficient to count the number of
rotationally distinct ways in which the temporal labels on the peripheral edges can be either greater or less than the
temporal label on the central edge.  

\begin{proposition}
\label{DiasterCountingProposition} Let $N=D(a,b)$ with central edge $c_{e}$ and peripheral edges $e_{1}, e_{2},
\ldots, e_{a+b}$. If $M\cong_{\mathrm{G}}N$ by $\phi:N \to M$, $\phi$ is a temporal isomorphism of $N$ and $M$ if and
only if for any $e_{i} \in N$ if $\tau(e_{i}) > \tau(c_{e})$, $\tau'(\phi(e_{i}))>\tau'(\phi(c_{e}))$, and if
$\tau(e_{i}) < \tau(c_{e})$, $\tau'(\phi(e_{i}))<\tau'(\phi(c_{e}))$
\end{proposition}
\emph{Proof} --- It is trivial to show that if  $\phi:N \to M$ is a temporal isomorphism, then for all $\tau(e_{i}) >
\tau(c_{e})$, it must be the case that $\tau'(\phi(e_{i}))>\tau'(\phi(c_{e}))$. If $\tau(e_{i}) >
\tau(c_{e})$, then $\langle c_{e}, e_{i} \rangle$ is a temporal path in $N$, and  $\langle \phi(c_{e}), \phi(e_{i})
\rangle$ is a temporal path in $M$. By Definition \ref{TemporalPathDefinition}, $\tau'(\phi(c_{e}))<\tau'(\phi(e_{i}))$.
Similarly whevener $\tau(e_{i}) < \tau(c_{e})$, it must be the case that $\tau'(\phi(e_{i}))<\tau'(\phi(c_{e}))$. 

To show the converse we assume that whenever $\tau(e_{i}) > \tau(c_{e})$, $\tau'(\phi(e_{i}))>\tau'(\phi(c_{e}))$.
Temporal paths in diasters can be at most three edges long: $P=\langle e_a, c_e, e_b \rangle$. If $\phi$ preserves
inequality relationships between $\tau(e_i)$ and $\tau(c_e)$, then it must be the case that $\tau'(\phi(e_a)) <
\tau'(\phi(c_a))<\tau'(\phi(e_b))$, so $\langle \phi(e_a), \phi(c_e), \phi(e_b) \rangle$ is a temporal path in M. Since
any temporal path in $N$ will be a subpath of $P$, it is clear that $\phi$ preserves all temporal paths. \qed
\\
\\
This Proposition places the isotemporal classes into one-to-one and onto correspondence with the classes generated by
considering two diasters equivalent if there exists an isomorophism that preserves the inequality relationships between
temporal labels on peripheral edges and the central edge. Thus to count isotemporal classes it is sufficient to count
the number of distinct inequality relationships.

\begin{theorem} 
\label{IsotemporalClassesOfDiasters} If $a+b>1$, the number of isotemporal classes of a diaster
$\mathcal{N}(D(a,b)))=ab+a+b+1$, and if $a=b$, $\mathcal{N}(D(a,a)))=\frac{1}{2} (a^2+3a+2)$.
\end{theorem}
\emph{Proof} --- Let us consider the case when $a < b$, and for convenience say that $a$ peripheral edges appear on the
``left'' while $b$ appear on the ``right.'' $T=\{1,2,\ldots,a+b+1\}$, and any of these labels could be applied to the
central edge. Given that temporal label $t$ is assigned to the central edge, there can be up to $t-1$ or $a$ (whichever
is fewer) labels on the left peripheral edges such that $\tau(e_i) < t$. Furthermore if we assume that $k$ left
peripheral edges have temporal labels less than $t$, this fixes the number right peripheral edges with temporal labels
less than $t$ at $t-1-k$. All other peripheral peripheral edges must have temporal labels that are greater than $t$. Any
diasters with $k$ left and $t-1-k$ right peripheral edges with temporal labels less than $t$ belong to the same
isotemporal class because $k!(a-k)!(t-1-k)!(b-t+1+k)!$ graphical isomorphisms preserving these $k$ left and $t-1-k$
right peripheral edges satisfy Proposition \ref{DiasterCountingProposition}.

If $t \leq a$, $0 \leq k<t$ left peripheral edges can have temporal labels less than $t$. If $a<t \leq b+1$, from $k=0$
to $k=a$ left peripheral edges can have temporal labels less than $t$. And if $b+1 < t \leq a+b+1$ (in this case more
than $b$ peripheral edges - i.e. more than 0 on the left - have to have temporal labels less than $t$), the number of
left peripheral edges with temporal labels less than $t$ must satisfy $a-(a+b+1-t) \leq k \leq a$. In each of these
three cases every unique $k$ corresponds to a distinct isotemporal class. Therefore:

$$
\mathcal{N}(D)=\sum_{t=1}^a ((t-1)+1) + \sum_{t=a+1}^{b+1} (a+1) +\sum_{t=b+2}^{a+b+1}(a-(a-(a+b+1-t))+1)$$
$$
=\frac{a(a+1)}{2}+ ((b+1)-(a+1)+1) (a+1) +\sum_{t=b+2}^{a+b+1}(a+b+2)-\sum_{t=b+2}^{a+b+1}t$$
$$
=\frac{a(a+1)}{2}+ (b-a+1) (a+1) +a(a+b+2)-(a(b+1)+\sum_{t=1}^{a}t)$$
$$
=\frac{a(a+1)}{2}+ (ab-a^2+a+b-a+1) +(a^2+ab+2a)-(ab+a)-\frac{a(a+1)}{2}$$
$$
\mathcal{N}(D(a,b))=ab+a+b+1$$

In considering the case when $a=b$ it is critical to note that the diaster is now symmetrical under reflection across
the central edge. Therefore networks with a total $t-1$ peripheral edges with temporal labels less than that on the
central edge, $k$ of which are on the left, are temporally isomorphic to networks with $t-1-k$ left peripheral edges
having the same property (by a graph isomorphism that maps left peripheral edges to the right and vice versa).
Consequently we need only to count the number of isotemporal classes with up to $\lfloor \frac{t-1}{2} \rfloor$ or $a$
(which ever is fewer) left peripheral edges with temporal labels less than that of the central edge.

If $1\leq t \leq a+1$, there is a distinct isotemporal class for all $0 \leq k \leq \lfloor \frac{t-1}{2} \rfloor $. If
$a+1 < t \leq 2a+1$, at least one left edge must have a temporal label less $t$. It is therefore sufficient to count the
number of ways to place temporal labels on the left peripheral edges with values \emph{greater} than $t$; it follows
that all other edges on the left must have temporal labels less than $t$. For any $t$ there will be $2a+1-t$ total edges
with temporal values greater than $t$. Assigning $0\leq k \leq \lfloor \frac{2a+1-t}{2} \rfloor$ to the left enumerates
all isotemporal classes when $t > a+1$. 

$$
\mathcal{N}(D) = \sum_{t=1}^{a+1} (\lfloor \frac{t-1}{2} \rfloor  +1) + \sum_{t=a+2}^{2a+1} (\lfloor
\frac{2a+1-t}{2} \rfloor 
+1) $$
$$
= \sum_{t=2,  t \mathrm{even}} ^{t \leq a+1} \frac{t}{2} +\sum_{t=1,  t \mathrm{odd}}^{t \leq a+1} \frac{t+1}{2} +
 \sum_{t=a+2,  t \mathrm{even}}^{t\leq 2a+1} \frac{2a-t+2}{2} + \sum_{t=a+2,  t \mathrm{odd}}^{t\leq
2a+1}\frac{2a-t+3}{2} $$
\\
\noindent Assume $a$ is even.
\\
$$
=(1+\cdots+a/2)+(1+\cdots+(a+2)/2)+(a/2+\cdots+1)+(a/2+\cdots+1)$$
$$
= \frac{(\frac{a}{2}+1)\frac{a}{2}}{2} +\frac{(\frac{a+2}{2}+1)\frac{a+2}{2}}{2} +
 \frac{(\frac{a}{2}+1)\frac{a}{2}}{2} +\frac{(\frac{a}{2}+1)\frac{a}{2}}{2} $$
$$
=\frac{1}{2}(3\frac{a+2}{2}\frac{a}{2} + \frac{a+4}{2}\frac{a+2}{2})$$
$$
=\frac{1}{2}\frac{(3a^2+6a) +(a^2+6a+8)}{4}$$
$$
\mathcal{N}(D(a,a))=\frac{1}{2}(a^2+3a+2)$$
\\
\noindent Assuming $a$ is odd,
\\
$$
\mathcal{N}(D)=(1+\cdots+\frac{a+1}{2})+(1+\cdots+\frac{a+1}{2})+(\frac{a-1}{2}+\cdots+1)+(\frac{a+1}{2%
}+\cdots+1)$$
$$
=\frac{1}{2}(3 \frac{a+1}{2}\frac{a+3}{2} + \frac{a-1}{2}\frac{a+1}{2})$$
$$
=\frac{1}{2}\frac{(3 a^2 + 12a + 9) + (a^2-1)}{4}$$
$$
\mathcal{N}(D(a,a))=\frac{1}{2}(a^2+3a+2) \qed$$
\\
The number of isotemporal classes of a diaster can be determined geometrically by a method reminiscent of calculating
${n \choose k}$ by finding the entry in the $n^{\mathrm{th}}$ row and $k^{\mathrm{th}}$ column of Pascal's Triangle.
Three summation terms enumerate $\mathcal{N}(D(a,b))$ in the proof of Theorem \ref{IsotemporalClassesOfDiasters}. The
indices of these terms span $t=1, 2, \ldots, a+b+1$. From $t=1, \ldots, a$, $t-1$ classes are identified for each $t$.
At each $a+1 \leq t \leq b+1$, $a+1$ more are added, and between $b+2 \leq t \leq a+b+1$, another $(a+b+1)+1-t$ are
added for each t. 

These components can be visualized on a square lattice where columns represent possible values of the temporal label on
the central edge ($t$), and rows represent the number of left peripheral edges with temporal labels less than $t$ ($k$),
starting at 0. Those values of $t$ and $k$ which correspond to possible temporal networks of $D(4,7)$ are given as an
example in Figure 3A. Generally, $\mathcal{N}(D(a,b))$ is equivalent to the number of square lattice vertices within a a
$45^\circ$ trapezoid with base length $a+b+1$ and height $a+1$. The area of such a trapezoid is
$\frac{1}{2}((a+b+1)+(a+b+1-2a))(a+1)=\frac{1}{2}(2ab+2b+2a+2)$.

The additional symmetry present in the $D(a,a)$ case reduces the number of isotemporal classes, and fewer distinct
isotemporal classes are enumerated for each possible value of $t$ (Figure 3B). $\mathcal{N}(D(a,a))$ can be easily
determined geometrically by noting that $\frac{1}{2}(a^2+3a+2) = \frac{(a+1)(a+2)}{2}=\sum_{i=1}^{a+1}i$.Thus
$\mathcal{N}(D(a,a))$ equals the number of vertices contained in a $45^\circ$ right triangle with legs of length $a+1$.

\section{Permutations preserving $\mathcal{N}(G)$} As we saw in Figure 1, it is often possible to permute temporal
labels and derive non-superimposable temporal networks that belong to the same isotemporal class. We can use this
observation to show that there exist mappings with the unusual property of preserving the number of isotemporal classes
while violating graphical isomorphism. 

\begin{proposition}
\label{LegalPermutations} Let $N\cong_{\mathrm{G}}M$ be temporal networks such that $N\not\cong_{\mathrm{L}}M$. If there
exists a permutation of temporal labels $\pi=p_1p_2\cdots p_n$ such that each $p_i$ is a sequential transposition of
temporal labels on non-adjacent edges with $\pi(N)\cong_{\mathrm{L}}M$, then $N\cong_{\mathrm{T}}M$ 
\end{proposition}
\emph{Proof} --- First we will show that if such a $\pi$ exists, $N\cong_{\mathrm{T}}M$. Let $\pi(N)\cong_{\mathrm{L}}M$
and $p_n=(a, a+1)$ which permutes the labels on edges $e_a=\tau^{-1}(a)$ and $e_{a+1}=\tau^{-1}(a+1)$. All temporal
paths including $e_a$ will be of the form $\langle e_1, \ldots, e_j, e_a, e_{j+1}, \ldots, e_{j+k} \rangle$ with
temporal labels $\langle t_1, \ldots, t_j, a, t_{j+1}, \ldots, t_{j+k} \rangle$ such that $t_i<a$ for all $i\leq j$ and
$t_i > a+1$ for all $i \geq j+1$. The value of $t_{j+1}$ must be greater than $a+1$ since the $e_{a+1}$ is not adjacent
to $e_a$.

Swapping the temporal labels on $e_a$ and $e_{a+1}$ means that the path $\langle e_1, \ldots,$ $e_j, e_a, e_{j+1},
\ldots, e_{j+k}\rangle$ will now have temporal labels $\langle t_1,
\ldots, t_j, a+1, t_{j+1}, \ldots,$ $t_{j+k} \rangle$. This is a temporal path since $a+1 < t_{j+1}$, and thus
$N\cong_{\mathrm{T}}p_n(N)$. By the same reasoning, $N\cong_{\mathrm{T}}p_{n-1}(p_n(N))$ and moreover,
$N\cong_{\mathrm{T}}\pi(N)$. Any isomorphism preserving edge labels will preserve temporal paths, so
$N\cong_{\mathrm{T}}M$. \qed
\\
\\
The converse of this Proposition is not generally true, however, in the case of diasters, it can be shown
with the help of the following Lemma. 

\begin{lemma}
\label{SequentialSwaps} Let $A=a_1, a_2, \ldots, a_n$ and $B=b_1, b_2, \ldots, b_n$ such that $a_i, b_i, \in \{0,1\}$
and  $|\{a_i\in A | a_i=0\}|=|\{b_i\in B | b_i=0\}|$. There exists a permutation on the elements of $A$
$\pi=p_1p_2\cdots p_n$ where $p_i$ permutes only non-identical sequential elements, such that $p(A)=B$.
\end{lemma}
\emph{Proof} --- We can construct the permutation $\pi$ as follows. Consider $a_1$ and $b_1$. If these are equal then
consider $a_2$ and $b_2$. If they are not equal, then let us assume without loss of generality that $a_1=0$ and $b_1=1$.
There must be be some $a_i=1$ for $i>1$ by the assumption that $A$ and $B$ contain the same number of 1s and 0s. Let us
say that $a_j=1$, and $j$ is the smallest value that satisfies this condition. The permutation
$p_1=(1,2)\cdots(j-2,j-1)(j-1,j)$ satisfies the proposition; initially $a_j=1$ and $a_{j-1}=0$ making $(j,j-1)$ a valid
permutation, applying $(j-2,j-1)$ is similarly valid because $a_{j-2}=0$ and $a_{j-1}$ now equals $1$. After the
application of all of $p_1$, $a_1=1$, and $a_i=0$ for $1<i\leq j$. We can now consider consecutive elements of $A$ and
$B$ starting with $a_2$ and $b_2$, each time applying the permutation $p_i=(k, k+1)\cdots(k+l-1,k+l)$ where $a_k$ is the
first element such that $a_k\neq b_k$ and $a_{k+l}$ is the first $a_i \in A$ such that $a_i=b_k$. When all $m$
differences between $A$ and $B$ have been corrected $p_m (\cdots (p_2 (p_1 (A))\cdots)=B$.\qed

\begin{proposition}
\label{DiasterSpecificPermutationConverse} Let $N\cong_{\mathrm{G}}M$ be temporal diasters such that
$N\cong_{\mathrm{T}}M$. There exists a permutation of temporal labels $\pi=p_1p_2\cdots p_n$ such that each $p_i$ is a
sequential transposition of temporal labels on non-adjacent edges and $\pi(N)\cong_{\mathrm{L}}M$.
\end{proposition}
\emph{Proof} --- $N$ and $M$ must have the same number of vertices to be temporally isomorphic, so let $N,M=D(l,m)$. For
convenience let us call those peripheral edges containing the vertex of degree $l+1$ the ``left edges.'' It follows from
Proposition \ref{DiasterCountingProposition} that  if  $N\cong_{\mathrm{T}}M$, they must have the same temporal value on
their central edges ($c_e$). Therefore the same number of edges in $N$ and $M$ will have temporal values less than
$\tau(c_e)$; let this number be denoted $g$. Let $A=a_1, a_2, \ldots , a_g$ and $B=b_1, b_2, \ldots , b_g$ such that
$a_i$ (or $b_i) = 0$ if the edge with temporal label $i$ in $N$ (or $M$) is found among the left edges. If the edge with
label $i$ is on the right, $a_i$ is assigned the value 1. $A$ and $B$ satisfy the conditions of Lemma
\ref{SequentialSwaps} implying the existence of a permutation $\pi_1=p_1p_2\cdots p_m$ such that $p(A)=B$, and each
$p_i$ is a product of permutations of non-identical sequential members of $A$.  

We need now to construct a similar permutation for those edges in $N$ and $M$ with temporal labels greater than
$\tau(c_e)$. This can be easily done by the method above. Let this permutation be denoted $\pi_2$, and permute the
representation ($A'$) of the position of those temporal labels in $N$ greater than $\tau(c_e)$ to the representation
($B'$) of the position of those same labels in $M$. 

The permutation $\pi=\pi_1 \pi_2$ satisfies this Proposition because the swapping non-identical of the elements in $A$
and $A'$ by each $p$ comprising $\pi$ corresponds to swapping a temporal label on a left edge for one on a right edge,
and therefore these edges are non-adjacent. Furthermore, given that $p$ swaps sequential elements $a_j$ and $a_{j+1}$,
when applied to temporal networks, it acts on sequential temporal labels. Finally $\pi(N)\cong_{\mathrm{L}}M$ because
$\pi_1$ places those temporal labels less than $\tau(c_e)$ in $N$ into the left edges (if they appear among the left
edges of $M$) and into the right if they appear among the right edges of $M$. Similarly $\pi_2$ permutes the temporal
labels greater than $\tau(c_e)$ in $N$ into the same arrangment as $M$. \qed
\\
\\
That two diasters $N$ and $M$ are temporally ismorophic if and only if there exists a permutation $\pi$ which is a
product of sequential transpositions on non-adjacent edges that takes the temporal labels of $N$ to $M$ means that all
temporal labelings that belong to a particular isotemporal class can be enumerated by exhaustively applying all possible
sequential transpositions on non-adjacent edges. 

We can take advantage of this fact to determine the number of isotemporal classes of a variety of other graphs that are
structurally related to diasters. 

\begin{theorem}
\label{NewGraphClasses} Let $N=\{V,E,T,\tau\}$ and $M=\{W, F ,T ,\tau' \}$ be temporal networks.
\begin{enumerate}
\item  If  any two networks $N$ and $N'$ of the form $\{V,E, T, \psi \}$ (where $\psi$ is any temporal labeling) have
the property that $N\cong_{\mathrm{T}} N'$ if and only if there exists a permutation $\pi$ composed of sequential
transpositions on non-adjacent edges, 
\item if there exists a mapping $\phi : V \to W$ such that $\phi(\{v_i,v_j\})=\{\phi(v_i),\phi(v_j)\} \in F$ if and only
if $\{v_i,v_j\}
\in E$,  $\phi(e_i)$ is adjacent to $\phi(e_j)$ in $M$ if and only if $e_i$ and $e_j$ are adjacent in $N$, and for some
$\tau(e_i)=\tau'(\phi(e_i))$, 
and
\item if $N$ and $\phi(N)$ have the same edge automorphisms.
\end{enumerate} then $\mathcal{N}(N)=\mathcal{N}(M)$.
\end{theorem}
\emph{Proof} --- Because $\phi$ is bijective on edges, $|E|=|F|$, and because these edges have the same automorphisms, 
$N$ and $M$ will have the the same number of symmetrically distinct possible temporal labelings. This number ($P_N$) can
be determined using the P\`{o}lya Enumeration Theorem. 

Because we can enumerate all $N'$ that are temporally isomorphic to a particular network $N$ by exhaustively permuting
all sequential temporal labels that appear on non-adjacent edges, we precisely determine what subset of the total $P_N$
labelings belong to this isotemporal class. Let $\pi$ be the particular permutation composed of sequential
transpositions of non-adjacent edges such that $\pi(N)\cong_{\mathrm{L}}N'$. Because $\phi$ preserves adjacencies and
labels of edges, $\pi$ is a isomorphism of $\phi(N)$. This is true for all $\pi$ and therefore the number of labelings
in the isotemporal class of $N$ is at least as large as the size of isotemporal class of $\phi(N)$. Because $N$ and $M$
have the same edge connectivity and symmetry, and $\phi$ preserves temporal labels, it cannot be the case that there
exists some $\pi$ and $\pi'$ such that $\pi(N)=N'$ and $\pi'(N)=N''$ where $N' \neq N''$, but
$\pi(\phi(N))=\pi'(\phi(N))$. Thus the isotemporal classes of $N$ and $\phi(N)$ have the same size. 

Since $P_N$ can be written as a sum of the sizes of the $\mathcal{N}(N)$ distinct isotemporal classes, $P_N=C_1+C_2
+\cdots + C_{\mathcal{N}(N)}$, and since $\phi$ preserves each $C_i$, the $P_N$ labelings of $M$ will be partitioned
into the same number of isotemporal classes as $N$. \qed
\\
\\
At least two classes of pseudograph satisfy Theorem \ref{NewGraphClasses} with respect to diasters. 

\begin{definition} --- A \textbf{beachball} $\beta(k)$ is a pseudograph comprising two vertices connected by $k$ edges.
\end{definition}

\begin{definition} --- A \textbf{daisy} $\delta(k)$ contains a single vertex and $k$ loops.
\end{definition}

Two daisies or beachballs can be joined together by a single edge connected to the single vertex of a daisy, or exactly
one of the vertices of the beacball to form a \textbf{stem-structure}, denoted by $S(G,H)$ where $G$ and $H$ are the
constituent graphs joined by the connecting edge. For example, the pseudograph with composed of two singly connected
vertices each with a loop is denoted $S(\delta(1),\delta(1))$. In this framework, a diaster $D(a,b)$ can be be viewed as
$S(\sigma(a),\sigma(b))$, where $\sigma(a)$ indicates a \textbf{star graph} with $a+1$ vertices. 

\begin{proposition}
\label{DaisiesAndBeachballs} For $a,b>0$, if $a\neq b$, $\mathcal{N}(S(G,H))=ab+a+b+1$ where $G=\sigma(a)$, $G=\beta(a)$
or $G=\delta(a)$, and $H=\sigma(b)$, $H=\beta(b)$ or $H=\delta(b)$; if $a=b$,
$\mathcal{N}(S(\sigma(a),\sigma(a)))=\mathcal{N}(S(\delta(a),\delta(a)))=\mathcal{N}(S(\beta(a),\beta(a)))=\frac{1}{2}
(a^2+3a+2)$
\end{proposition}
\emph{Proof} --- This follows directly from the formula for $\mathcal{N}(D(a,b))$ (Theorem
\ref{IsotemporalClassesOfDiasters}), and Theorem
\ref{NewGraphClasses} if we let $N$ be a diaster and $M$ the new class of graph (or pseudograph). Inspection of the
properties of these new classes of pseudograph shows that the conditions of the latter Theorem are satisfied: First,
as in the case of diasters, the fact that $S(G,H)$ or $S(G,G)$ contain a central edge through which establishes two
mutually non-adjacent sets of left and right edges means that Proposition \ref{DiasterCountingProposition} holds, and
Proposition \ref{DiasterSpecificPermutationConverse} follows. 
Therefore, two labelings of  $S(G,H)$ or $S(G,G)$ are temporally isomorphic if and only if there exists a permutation
composed of sequential transpositions on non-adjacent edges between them edge.

Secondly, any function that takes the $a$ left peripheral edges of $D(a,b)$ to the edges of $G$, the remaining $b$
peripheral edges to the edges of $H$, the central edge of $D(a,b)$ to that of $S(G,H)$, and preserves temporal labels
will satisfy the requirements of $\phi$. 

Lastly, $D(a,b)$ and $S(G,H)$ both have edge automorphism groups isomorphic to $S_a \times S_b$ where $S_n$ is the
symmetric group on $n$ objects. $D(a,a)$ and $S(G,G)$ both have $S_a$ symmetry among both their left edges and right
edges, as well as reflective symmetry across their central edges. If we label the left edges as $\{e_1, e_2, \ldots,
e_a\}$ and the right edges as $\{e_{a+1}, e_{a+2}, \ldots,e_{2a}\}$, the edge automorphism group of  $D(a,a)$ (and
$S(G,G)$) is generated by $\{(e_1,e_2),(e_2,e_3), \ldots, (e_{a-1},e_a),(e_a,e_1)\} \cup
\{(e_{a+1},e_{a+2}),(e_{a+2},e_{a+3}),
\ldots,$ $(e_{2a-1},$ $e_{2a}),(e_{2a},e_{a+1})\} \cup \{(e_1,e_{a+1})(e_2,e_{a+2})\cdots(e_a,e_{2a})\}$, and is
isomorphic to $S_a \times S_a \times \mathbb{Z}_2$. \qed
\\
\\
This result greatly expands the number of graphs with closed formulas for their number of isotemporal classes (Figure
4). When the number of left and right edges is unequal there are six distinct classes: diasters, dibeachballs,
didaisies, aster-dasies, aster-beachballs and daisy-beachballs. When the number of left and right edges are identical,
this result only implies the number of isotemporal classes of diasters, dibeachballs, didaisies, since in the other
classes of pseudographs there is a lack of reflective symmetry across the central edge. 

Theorem 2 is particularly useful considering that there is no general method to determine closed formulae for the number
of isotemporal classes of a particular type of graph. To be able to determine $\mathcal{N}(G)$ for new types of graph
simply by noting that 1) any two networks are in the same isotemporal class if and only if there exists a permutation
composed of transpositions of sequential edge labels taking one network to the other, 2) the two types of graph have the
same edge adjacencies, and 3) the two types of graph have the same edge automorphisms --- is a significant step toward
generalized approaches to determining $\mathcal{N}(G)$ for arbitrary graphs.

\begin{figure}
  \begin {center}
\resizebox{0.6\textwidth}{!}{%
  \includegraphics{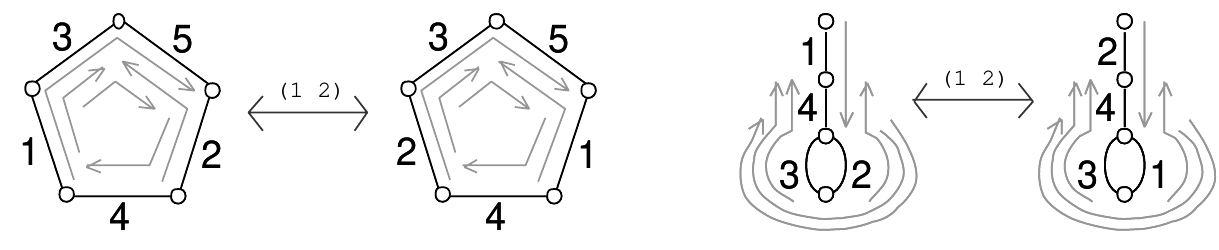}
}
  \end {center}
\caption{\emph{Examples of temporal isomorphism} --- The two 5-cycles on the left have identical temporal paths (gray)
but
different temporal labels, as do the multigraphs at right. Both pairs of network are temporally isomorphic by the edge
label permutation $(1, 2)$.}
\label{fig1}       
\end{figure}

\begin{figure}
\begin{center}
\resizebox{0.5\textwidth}{!}{%
  \includegraphics{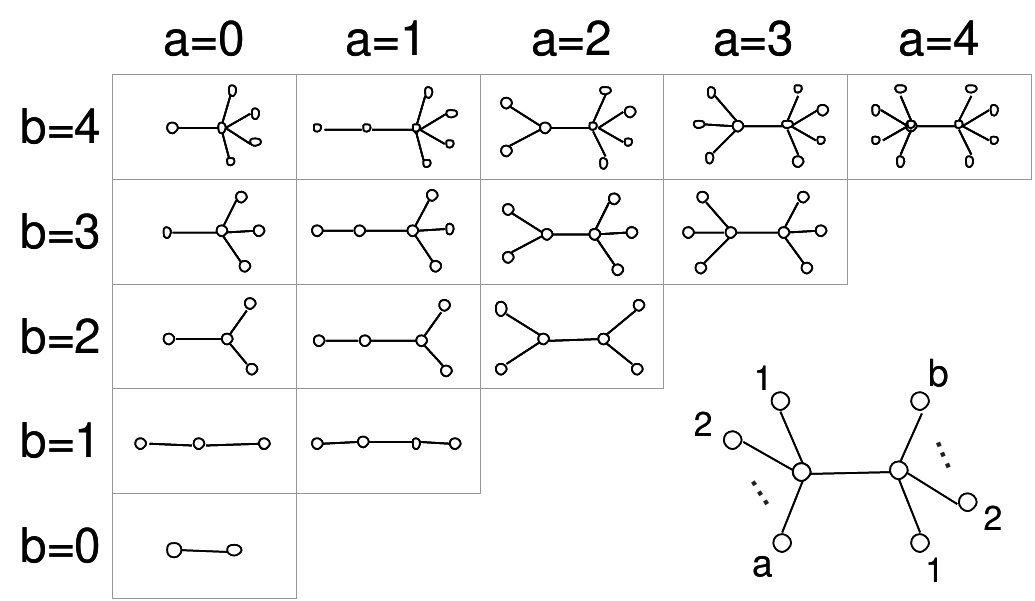}
}
  \end {center}
\caption{\emph{The class of diaster graphs} --- Diasters are shown for varying numbers of right and left peripheral edges
(a and b).}
\label{fig2}       
\end{figure}

\begin{figure}
\begin{center}
\resizebox{0.6\textwidth}{!}{%
  \includegraphics{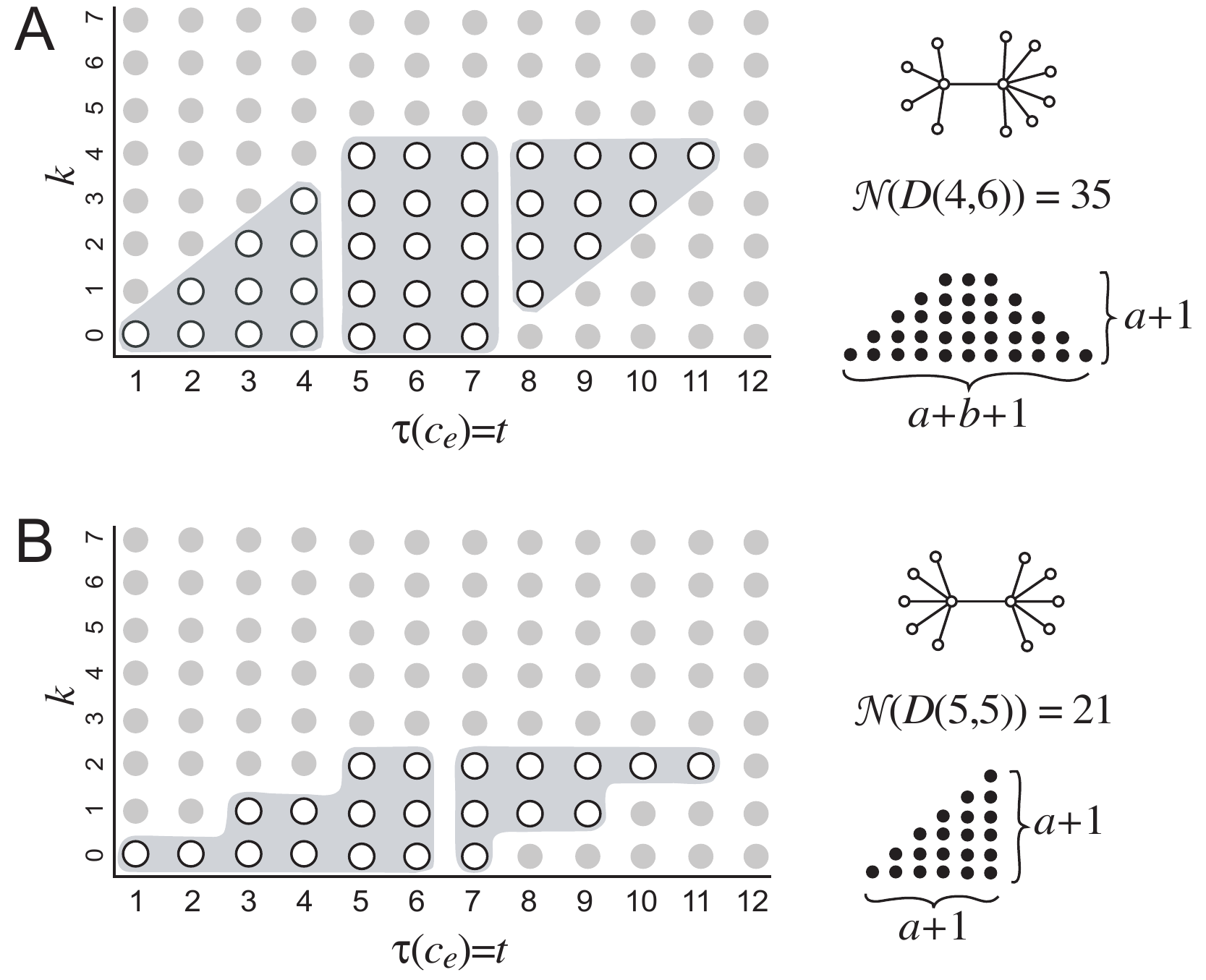}
}
  \end {center}
\caption{\emph{ Geometric interpretation of the number of isotemporal classes of
diasters $D(a,b)$} --- If $a\neq b$, $\mathcal{N}D(a,b) = ab + a + b + 1$, if $a=b$, $\mathcal{N}D(a,b) =
\frac{1}{2}(a^2+3a+2)$. These values can be determined by counting the number of nodes in a square lattice within
45$^\circ$ and 22.5$^\circ$ parallelograms.}
\label{fig3}       
\end{figure}

\begin{figure}
\begin{center}
\resizebox{0.8\textwidth}{!}{%
  \includegraphics{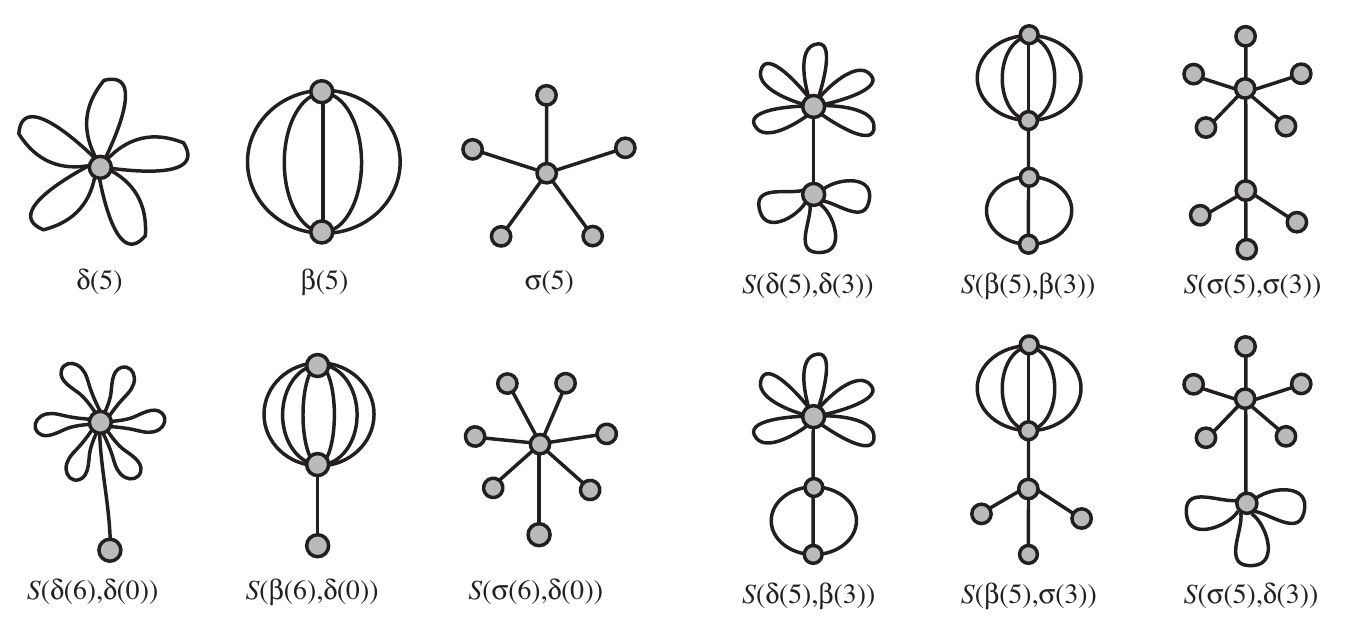}
}
\end{center}
\caption{\emph{Pseduographs whose number of isotemporal classes follow from the diaster fomula and Theorem 2} ---
Clockwise from top left, daisies, beachballs and asters have a single isotemporal class; didaisies, diasters,
dibeachballs, stemmed daisies, stemmed beachballs, stemmed asters (also asters), daisy-beachballs, aster-beachballs,
and aster-daisies all obey the formulae for the number of isotemporal classes of diasters.}
\label{fig4}       
\end{figure}

\end{document}